%% file: ak4.tex
\documentclass{article}
\title{Algebraic reduction of certain almost K\"ahler manifolds
}
\author{Paul-Andi Nagy}
\date{\today}
\usepackage{amsfonts,amssymb}
\usepackage{minitoc_href}

\newtheorem{teo}{Theorem}[section]
\newtheorem{lema}{Lemma}[section]
\newtheorem{pro}{Proposition}[section]
\newtheorem{defi}{Definition}[section]
\newtheorem{rema}{Remark}[section]
\newtheorem{coro}{Corollary}[section]
\newtheorem{nr}{}[section]
\begin{document}
\maketitle
\abstract{\normalsize We study almost K\"ahler manifolds whose curvature tensor satisfies the third
curvature condition of Gray. We show that the study of manifolds within this class reduces to the study 
of a subclass having the property that the torsion of the first canonical Hermitian connection 
has the simplest possible algebraic form. 
This allows to understand the structure of the K\"ahler nullity of an almost K\"ahler manifold with parallel torsion.
\large
\tableofcontents
\section{Introduction}
An  almost K\"ahler manifold (shortly ${\cal{AK}}$) is a Riemannian manifold $(M^{2n},g)$, together with a compatible almost complex structure 
$J$, such that the K\"ahler form $\omega=g(J \cdot, \cdot)$ is closed. Hence, almost K\"ahler geometry is nothing else that symplectic 
geometry with a prefered metric and almost complex structure. Since symplectic manifolds often arise in this way is rather 
natural to ask under which conditions on the metric we get integrability of the almost complex $J$. In this direction, 
a famous (still open) conjecture of S. I. Golberg asserts that every compact, Einstein, almost K\"ahler 
manifold is, in fact, K\"ahler. Nevertheless, they are a certain number of partial results, supporting this 
conjecture. First of all, K. Sekigawa proved \cite{Seki1} that 
the Goldberg conjecture is true when the scalar curvature is positive.
We have to note that the Golberg conjecture is definitively not true with the compacity assumption removed 
by the examples of \cite{Apo2} in complex dimension $n \ge 3$ and those of 
\cite{Apo, arm2,Nur} in (real) dimension $4$.  
In \cite{Apo2}, a potential source of compact almost K\"ahler, Einstein manifolds is considered, 
namely those compact K\"ahler manifolds whose Ricci tensor admits two distinct, constant eigenvalues; integrability 
is proven under certain positivity conditions. 
The rest of known results, most of them 
enforcing or replacing the Einstein condition with some other natural curvature assumption are mainly 
in dimension $4$. To cite only a few of them, we mention the beautiful series of papers \cite{Apo1, Apo3, Apo4} giving 
a complete local and global classfication of almost K\"ahler manifolds of $4$ dimensions satisfying the second and third 
Gray condition on the Riemannian curvature tensor. Other recent results, again in $4$-dimensions, are concerned with the study of local obstructions 
to the existence of Einstein metrics \cite{arm2}, $\star$-Einstein metrics \cite{Seki2}, etc. As for the class of 
almost K\"ahler manifolds satisfying the second curvature condition of Gray it was established recently 
\cite{Nagy, Nagy3} that the it coincides with the class of almost K\"ahler manifolds whom torsion of 
the first canonical connection is parallel. In view of the known twistorial examples 
\cite{Davidov, Ivanov}, we have therefore a strong ressemblence between this class and that of nearly K\"ahler 
manifolds. 
\par
In this paper our main object of study will be the class of almost K\"ahler manifolds satisfying the third curvature 
condition of Gray (shortly ${\cal{AK}}_3$). 
To the best of our knowledge, in dimension greater than $6$ almost nothing is known about the structure 
of this class of manifolds.
Our approach to the study of the class ${\cal{AK}}_3$-manifolds will be directed from the point 
of view of the canonical Hermitian 
connection. Actually, we are going to study the geometric as well as the algebraic effects of the third curvature 
condition of Gray over the torsion of the last mentioned connection. 
\begin{teo}
Let $(M^{2n},g,J)$ belong to the class ${\cal{AK}}_3$. Then, over some dense open subset, the manifold $M$ is locally the 
Riemannian product of a strict 
almost K\"ahler manifold with parallel torsion and a special ${\cal{AK}}_3$-manifold. 
\end{teo}
The precise definition of special ${\cal{AK}}_3$ is given in section 4. They are caracterized by an algebraic feature of 
their K\"ahler nullity which generalizes, to some extent, properties of almost K\"ahler $4$-manifolds. Finally, 
theorem 1.1 enables us to give a first structure result concerning almost K\"ahler manifolds with parallel torsion. 
\begin{coro}
Let $(M^{2n},g,J)$ be almost K\"ahler with parallel torsion. Then $M$ is locally the Riemannian product of a 
K\"ahler manifold, a strict ${\cal{AK}}$-manifold with parallel torsion and 
a locally $3$-symmetric space of type $I$.
\end{coro}
Here, and in view of further progress on the classification of almost K\"ahler manifolds with parallel torsion, 
we regrouped in the class I those special almost K\"ahler manifolds with parallel torsion. \par
The paper is organised as follows. In section 2 we review some classical facts and definitions 
from almost Hermitian geometry, mainly related to the first canonical Hermitian connection and 
Gray's curvature conditions. The main differential geometry properties of the torsion of the canonical connnection, 
such as the parallelism in directions orthogonal to the K\"ahler nullity are established in section 3. The last section 
is devoted to the proof of theorem 1.1. We start from a result used in the ${\cal{AK}}_2$-context (see \cite{Nagy}) 
and observe that it continues to hold in the present setting. Then, due the to more flexibility of 
${\cal{AK}}_3$ we need a finer analysis in order to produce the proof of theorem 1.1. One of main ingredients 
consists in using a suitablly defined Hermitian Ricci tensor in order to eliminate those parts of the 
tangent bundle of the manifold obstructing the decomposition in  theorem 1.1. 
\section{Preliminaries}
Let $(M^{2n},g,J)$ be an almost Hermitian manifold. Let $\nabla$ be the Levi-Civita connection of the metric $g$. Then for all $X$ in $TM$ we have that 
$\nabla_XJ$ is a skew-symmetric, $J$-anti-invariant endomorphism of $TM$. Imposing further algebraic constraints to the tensor $\nabla J$ leads to 
the introduction of several classes of almost Hermitian manifolds. We will recall the definition of those that will be significant for this paper.\par
For instance, $(M^{2n},g,J)$ is quasi-K\"ahler iff for all $X,Y$ in $TM$ we have : 
$$ (\nabla_{JX}J)JY=-(\nabla_XJ)Y.$$
If $\omega=g(J \cdot, \cdot)$ is the K\"ahler form of $(g,J)$ then we are in presence of an almost K\"ahler structure iff $d \omega=0$. It should 
be noted that almost K\"ahler structures are always quasi-K\"ahler. \par
An important object related to the almost Hermitian structure $(g,J)$ is the first canonical Hermitian connection defined 
by 
$$ \overline{\nabla}_XY=\nabla_XY+\eta_XY$$
whenever $X,Y$ are vector fields on $M$ where, 
to save space, we setted $\eta_XY=\frac{1}{2}(\nabla_XJ)JY$. We obtained a metric Hermitian connection on 
$M$, that is $\overline{\nabla}g=0$ and $\overline{\nabla}J=0$. 
The torsion tensor of the canonical Hermitian canonical connection, to be denoted by $T$ is given by 
$$T_XY=\eta_XY-\eta_YX$$
for all $X,Y$ in $TM$. Then by the {\bf{torsion}} of the almost Hermitian manifold $(M^{2n},g,J)$ we will mean 
simply the torsion tensor of the canonical Hermitian connection. Note that among connections respecting both 
the metric $g$ and the almost complex structure $J$ the first canonical Hermitian connection is the one 
whose torsion has minimal possible norm \cite{Gaud}.
\par
For the almost Hermitian $(M^{2n},g,J)$ be almost K\"ahler one requires that the K\"ahler form 
$\omega(X,Y)=<JX, Y>$ be closed. In our notations, this is equivalent to have 
\begin{nr} \hfill
$<T_XY,Z>=-<\eta_ZX,Y> \hfill $
\end{nr}
for all $X,Y,Z$ in $TM$. When dealing with almost K\"ahler manifolds, this relation will be used almost implicitely in  
this paper. \par  
We denote by $R$ resp. $\overline{R}$ the curvature tensor of the Levi-Civita connection $\nabla$ resp. of the 
canonical Hermitian connection $\overline{\nabla}$ with the convention that
$\overline{R}(X,Y)=-[\overline{\nabla}_X, 
\overline{\nabla}_Y]+\overline{\nabla}_{[X,Y]}$ for all vector fields $X,Y$ on $M$.
Now a standard calculation involving the definitions yields to 
\begin{nr} \hfill 
$ \overline{R}(X,Y)Z=R(X,Y)Z+[\eta_X,\eta_Y]Z-\biggr [ d_{\overline{\nabla}}u(X,Y) \biggl ] Z \hfill $
\end{nr}
where 
$$ \biggl [ d_{\overline{\nabla}}u(X,Y) \biggr ] Z =(\overline{\nabla}_X\eta)(Y,Z)-(\overline{\nabla}_Y\eta)(X,Z)+
\eta_{T_XY}Z$$
for all vector fields $X,Y,Z$ on $M$. Note that $d_{\overline{\nabla}}u(X,Y)$ is a $J$-anticommuting endomorphism 
of $TM$, whenever $X,Y$ are tangent vectors to $M$.
\begin{rema}
In the formula (2.2), the notation $u$ stands for the tensor $\eta$, considered as a $1$-form with values in the 
bundle $\Omega^2(M)$. Then $d_{\overline{\nabla}}$, with its expression given below, is the twisted differential 
acting on twisted one forms, when considering the tangent bundle of $M$ endowed with the connection 
$\overline{\nabla}$. Since our discussion is intended to be self-contained and at the elementary level, we will 
keep things at the level of the notation after (2.2).
\end{rema}
The fact that $\overline{\nabla}$ is a Hermitian connection implies that $\overline{R}(X,Y,JZ,JU)=\overline{R}(X,Y,Z,U)$. Using this in formula (2.2), together with 
the skew-symmetry of the $J$-anticommuting endomorphism $\eta_X$ gives us : 
\begin{nr} \hfill 
$ \begin{array}{lr}
R(X,Y,Z,U)-R(X,Y,JZ,JU)=2\biggl [ (d_{\overline{\nabla}}u)(X,Y) \biggr ](Z,U).
\end{array}
\hfill $
\end{nr}
Using the symmetry property of the Riemannian curvature we also deduce that 
\begin{nr} \hfill 
$ \begin{array}{lr}
R(X,Y,Z,U)-R(JX,JY,Z, U)=2\biggl [ (d_{\overline{\nabla}}u)(Z,U) \biggr ](X,Y)
\end{array}
\hfill $
\end{nr}
whenever $X,Y,Z,U$ are in $TM$. We are going to make an intensive use of the two formulas below in the next section.
\par
The rest of the present section is destinated to present some classes of almost Hermitian manifolds, the third of which will be the object 
of our attention in this paper. 
We begin by recalling how one 
can distinguish several classes of almost Hermitian manifolds by "the degree of ressemblance" of their 
Riemannian curvature tensor with the curvature tensor of a K\"ahler manifold \cite{Gray1, Tri2}: \\
$\begin{array}{lr} 
(G_1) : \ R(X,Y,JZ,JU)=R(X,Y,Z,U) \\
(G_2) : \  R(X,Y,Z,U)-R(JX,JY,Z,U)=R(JX,Y,JZ,U)+R(JX,Y,Z,JU) \\
(G_3) : \ R(JX,JY,JZ,JU)=R(X,Y,Z,U)
\end{array} $
$\\$
Using the first Bianchi identity it is a simple exercise to see that $G_1 \Rightarrow G_2 \Rightarrow G_3$. It is also clear that a K\"ahler 
structure satisfies all the three conditions.  Let us set now some notations. \par
Following \cite{Gray1}, let ${\cal{AK}}$ be the class of almost K\"ahler manifolds. Then the class ${\cal{AK}}_i, 1 \le i \le 3$ contains those 
almost K\"ahler manifolds whose curvature tensor satisfies the condition $(G_i)$. Obviously, we have the inclusions : 
$$ {\cal{AK}}_1 \subseteq {\cal{AK}}_2 \subseteq {\cal{AK}}_3. $$
Note that it was shown in \cite{Golberg} that locally $ {\cal{AK}}_1={\cal{K}}$, where ${\cal{K}}$ denotes the class of K\"ahler manifolds. Let us give 
a very simple proof of this result.
\begin{pro}Any almost K\"ahler manifold $(M^{2n},g,J)$ satisfying condition $(G_1)$ is K\"ahler.
\end{pro}
{\bf{Proof}} : \\
From (2.3) we obtain that condition $(G_1)$ is equivalent with 
$$ (\overline{\nabla}_X\eta)(Y,Z)- (\overline{\nabla}_Y\eta)(X,Z)+\eta_{T_XY}Z=0$$ 
for all $X,Y,Z$ in $TM$. Using the quasi-K\"ahler condition, one notices that the second term of this equation is $J$-anti-invariant in $X,Z$ whilst 
the last term is $J$-invariant in the same variables. Therefore we obtain 
$$ (\overline{\nabla}_{JX}\eta)(Y,JZ)+(\overline{\nabla}_X\eta)(Y,Z)=-2\eta_{T_XY}Z.$$
Using again the quasi-K\"ahler condition we find that left  hand side of the equation below is 
$J$-invariant in $X,Y$ whilst the right hand one is $J$-anti-invariant in the same variables. We conclude that $\eta_{T_XY}Z=0$ for all $X,Y,Z$ in 
$TM$ and then the almost K\"ahler condition (2.1) implies the vanishing of $T$, and hence that of $\eta$.
$\blacksquare$ \\ \par
The other inclusions between the previously defined classes are strict in dimensions $2n \ge 6$, as shows the examples of 
\cite{Davidov}, multiplied by K\"ahler manifolds. In the same spirit, the class ${\cal{AH}}_i, 1 \le i \le 3$ contains 
those almost Hermitian manifolds whose Riemannian curvature tensor satisfies condition $(G_i)$. By contrast 
with the almost K\"ahler case, note that they are classes of almost complex manifolds where conditions 
$(G_2)$ and $(G_3)$ are known to be equivalent, such as the class of Hermitian manifolds \cite{Gray1} and the 
class of locally conformally K\"ahler manifolds \cite{Ornea}.
\\ \par
To give a very simple interpretation of the third Gray condition on curvature, let us recall 
that we have a decomposition of the bundle of (real valued ) two forms : 
\begin{nr} \hfill
$ \Lambda^2(M)=\Lambda^{1,1}(M) \oplus [[ \Lambda^{2,0}(M)]] \hfill $
\end{nr}
into $J$-invariant and $J$-anti-invariant parts, where the action of $J$ on a two form $\alpha$ is given by 
$(J\alpha)(X,Y)=\alpha(JX,JY)$ for all $X,Y$ in $TM$.
It is easily checked that $(g,J)$ belongs to the class ${\cal{AH}}_3$ iff 
the Riemannian curvature operator of $g$ preserves the decomposition (2.5). \\
\section{The third Gray condition}
\input{gthree.tex}
\section{An algebraic decomposition}
\input{algak4f.tex}

\begin{flushright}
Paul-Andi Nagy \\
Institut de Math\'ematiques \\
rue E. Argand 11, CH-2007, Neuch\^atel \\ 
email : Paul.Nagy@unine.ch
\end{flushright}
\end{document}

%% file: gthree.tex
This section is dedicated to give an interpretation, in terms of the torsion of the first canonical Hermitian connection of the third Gray condition 
on curvature. In the context of quasi-K\"ahler geometry we will deduce some important geometric consequences of the condition 
$(G_3)$ putting the basis of the dicussion in the next section, where we will restrict finally to our main object of interest, the subclass 
of almost K\"ahler manifolds. We begin by the following 
very simple observation.
\begin{lema}The almost Hermitian manifold $(M^{2n},g,J)$ satisfies the condition $(G_3)$ iff 
\begin{nr} \hfill 
$ \biggl [ (d_{\overline{\nabla}}u)(Z,U) \biggr ](X,Y)=
\biggl [ (d_{\overline{\nabla}}u)(X,Y) \biggr ](Z,U).\hfill $
\end{nr}
Moreover, we have $(d_{\overline{\nabla}}u)(JX,JY)+(d_{\overline{\nabla}}u)(X,Y)=0$, hence 
$d_{\overline{\nabla}}u$ defines a symmetric endomorphism of $\Lambda_{-}^2(M)$.
\end{lema}
{\bf{Proof}} : \\
We change $Z$ and $U$ in $JZ$ and $JU$ respectively in (2.4) and take the sum with (2.3). Using $(G_3)$, 
we get  $ \biggl [ (d_{\overline{\nabla}}u)(JZ, JU) \biggr ](X,Y)=-
\biggl [ (d_{\overline{\nabla}}u)(X,Y) \biggr ](Z,U)$. To get (3.1) we change again $Z$ in $JZ$ and $U$ in $JU$ and 
use that $du_{\overline{\nabla}}(X,Y)$ is $J$-anticommuting. The rest is straighforward.
$\blacksquare$ \\ \par
The previous lemma enable us to give a number of useful algebraic properties of 
the tensor $\overline{R}$. 
\begin{coro}
Let $(M^{2n},g,J)$ be a quasi-K\"ahler manifold in the class ${\cal{AH}}_3$ and let 
$X,Y,Z,U$ be vector fields on $M$. The following holds : \\
(i) 
$$ \overline{R}(X,Y,Z,U)-\overline{R}(Z, U, X, Y)=<[\eta_X, \eta_Y]Z, U>-<[\eta_Z, \eta_U]X,Y>.$$
(ii) $\overline{R}(JX,JY),Z,U)=\overline{R}(X,Y,Z,U)$. \\
\end{coro}
{\bf{Proof}} : \\
Property (i) follows immediately from (3.1), the symmetry of Riemannian curvature operator and lemma 3.1. 
Since $\overline{R}$ is a Hermitian connection we have $\overline{R}(Z,U,JX,JY)=\overline{R}(Z,U,X,Y)$, hence 
(ii) follows from (i) and the quasi-K\"ahler condition on the tensor $\eta$.
$\blacksquare$ \\ \par
Note that the previous corollary is well known for nearly-K\"ahler manifolds (see \cite{Gray2} for instance). Also 
note that property (i) holds in fact for any almost Hermitian manifold in the class ${\cal{AH}}_3$. \par
With this preliminaries in mind, we are now going to investigate some illuminating consequences of the third Gray 
condition.
\begin{pro}
Let $(M^{2n}, g, J)$ be a quasi-K\"ahler manifold. If $(M^{2n},g,J)$ satisfies condition $(G_3)$ then : 
\begin{nr} \hfill
$ (\overline{\nabla}_{JX} \eta)(JY,Z)+(\overline{\nabla}_X \eta)(Y,Z)=0.\hfill $
\end{nr}
and 
$$ <\eta_{T_XY}Z,W>=<\eta_{T_ZW}X,Y>$$
whenever $X,Y,Z$ and $U$ are in $TM$.
\end{pro}
{\bf{Proof}} : \\
From lemma 3.1 we know that 
$$\begin{array}{lr}
<(\overline{\nabla}_X \eta)(Y,Z)-(\overline{\nabla}_Y\eta)(X,Z),U>+<\eta_{T_XY}Z,U>=\\
<(\overline{\nabla}_Z \eta)(U,X)-(\overline{\nabla}_U\eta)(Z, X),Y>+<\eta_{T_ZU}X,Y>
\end{array}$$
for all $X,Y,Z,U$ in $TM$. The quasi-K\"ahler condition ensures that the second term of each side is $J$-anti-invariant 
in $X$ and $Z$. Therefore, identifying the invariant parts (in $X$ and $Z$) of the previous equation we get 
$$\begin{array}{lr}
<(\overline{\nabla}_{JX} \eta)(Y,JZ)+(\overline{\nabla}_{X} \eta)(Y,Z),U>+2<\eta_{T_XY}Z,U>=\\
<(\overline{\nabla}_{JZ} \eta)(U,JX)+(\overline{\nabla}_{Z} \eta)(U,X),Y>+2<\eta_{T_ZU}X,Y>.
\end{array}$$
After rearanging terms and using the quasi-K\"ahler condition we find that 
$$\begin{array}{lr}
<(\overline{\nabla}_{JX} \eta)(JY,Z)+(\overline{\nabla}_{X} \eta)(Y,Z),U>=\\
<(\overline{\nabla}_{JZ} \eta)(JU,X)+(\overline{\nabla}_{Z} \eta)(U,X),Y>+2<\eta_{T_ZU}X,Y>-2<\eta_{T_XY}Z,U>.
\end{array}$$
But the left hand side of this equation is $J$-antinvariant in $Z$ and $U$ whilst the right hand is $J$-invariant (here 
one uses the quasi-K\"ahler condition on the last two terms) in the same variables. This leads to the proof of the 
proposition.
$\blacksquare$ 
\begin{rema}
The algebraic constraint in proposition 3.1 is automatically satisfied in both nearly K\"ahler and 
almost K\"ahler cases. 
\end{rema}
We are now going to show that the torsion of ${\cal{AH}}_3$-manifolds enjoys a particularly pleasant 
property.
\begin{pro}
Let $(M^{2n},g,J)$ be quasi-K\"ahler in the class ${\cal{AH}}_3$. Then : 
\begin{nr} \hfill 
$ \overline{\nabla}_{T_XY} \eta=0\hfill $
\end{nr}
for all $X,Y$ in $TM$. 
\end{pro}
{\bf{Proof}} : \\
Derivating (3.2) we obtain : 
$$ (\overline{\nabla}_{X,JY}^2 \eta)(JZ, \cdot)+(\overline{\nabla}^2_{X,Y}\eta)(Z, \cdot)=0.$$
Replacing $X$ with $JY$ and $Y$ with $JX$ respectively we obtain also 
$$ -(\overline{\nabla}_{JY,X}^2 \eta)(JZ, \cdot)+(\overline{\nabla}^2_{JY,JX}\eta)(Z, \cdot)=0.$$
After addition of these two equations and by making use of the Ricci identity (with respect of the connection 
$\overline{\nabla}$) we obtain : 
\begin{nr} \hfill 
$ \biggl [ \overline{R}(X,JY).\eta \biggr ](JZ, \cdot)+\overline{\nabla}_{T_X(JY)} \eta(JZ, \cdot)=
 (\overline{\nabla}_{X,Y}^2 \eta)(Z, \cdot)+(\overline{\nabla}^2_{JY,JX}\eta)(Z, \cdot).\hfill $
\end{nr}
Here, the action $G.\eta$ of an endomorphism $G$ of $TM$ on the tensor $\eta$ is defined by $(G.\eta)(X,Y)=G(\eta_XY)-\eta_{GX}Y-\eta_X(GY)$ 
for all $X,Y$ in $TM$. \par
We antisymmetrize (3.4) in $X$ and $Y$ in order to get, after using twice the Ricci identity : 
$$ \begin{array}{lr}
\biggl [ \overline{R}(X,JY)+\overline{R}(JX,Y). \eta \biggr ] (JZ, \cdot)+
(\overline{\nabla}_{T_X(JY)+T_{JX}Y} \eta)(JZ, \cdot)=\\
\biggl [ \overline{R}(JX,JY)-\overline{R}(X,Y). \eta \biggr ] (Z, \cdot)+
(\overline{\nabla}_{T_{JX}(JY)-T_{X}Y} \eta)(Z, \cdot).
\end{array} $$
Using corollary 3.1, (ii) it follows that $\overline{\nabla}_{T_X(JY)+T_{JX}Y} \eta)(JZ, \cdot)=
(\overline{\nabla}_{T_{JX}(JY)-T_{X}Y} \eta)(Z, \cdot)$ or further 
$(\overline{\nabla}_{JT_XY} \eta)(JZ, \cdot)=(\overline{\nabla}_{T_{X}Y} \eta)(Z, \cdot)$ since $(M^{2n},g,J)$ is 
quasi-K\"ahler. Our claimed result is implied now by (3.2).
$\blacksquare$\\ \par
We finish this section by showing that in the quasi-K\"ahler, ${\cal{AH}}_3$-setting, one can get a simplified form of 
the first Bianchi identity for the first canonical Hermitian connection . 
\begin{lema}
For any quasi-K\"ahler, ${\cal{AH}}_3$-manifold $(M^{2n},g,J)$ the following holds : 
$$ \sigma_{X,Y,Z}\biggl [ \overline{R}(X,Y)Z+T_{T_XY}Z \biggr ]=0 $$
whenever $X,Y,Z$ are in $TM$.
\end{lema}
{\bf{Proof}} : \\
We know that 
$$ <(\overline{\nabla}_X \eta)(Y,Z)-(\overline{\nabla}_Y \eta)(X,Z),U>=<(\overline{\nabla}_Z \eta)(U,X)-(\overline{\nabla}_U \eta)(Z,X),Y>.$$
Using the almost K\"ahler condition $<\eta_UX,Y>=-<T_XY, U>$ in the right hand side of the previous equation we are lead to 
$$<(\overline{\nabla}_X \eta)(Y,Z)-(\overline{\nabla}_Y \eta)(X,Z),U>=-<(\overline{\nabla}_Z T)(X,Y)-(\overline{\nabla}_U \eta)(Z,X),Y>.
$$
Taking the symmetric sum of the last identity yields to 
$$\begin{array}{lr}
2\sigma_{X,Y,Z} <(\overline{\nabla}_XT)(Y,Z),U>=\\
-\biggl [<(\overline{\nabla}_U \eta)(Z,X), Y>+
<(\overline{\nabla}_U \eta)(X,Y), Z>+<(\overline{\nabla}_U \eta)(Y,Z), X>
\biggr ].\end{array}
$$
But $<(\overline{\nabla}_U \eta)(Z,X), Y>+
<(\overline{\nabla}_U \eta)(X,Y), Z>=<(\overline{\nabla}_UT)(Z,X), Y>=-<(\overline{\nabla}_U \eta)(Y,Z), X>$ by the use of the almost K\"ahler condition 
$<T_ZX,Y>+<\eta_YZ,X>$ and using the first Bianchi identity for the 
connection $\overline{\nabla}$ finishes the proof of the lemma.
$\blacksquare$

%% file: algak4f.tex
In this section we begin the geometric study of ${\cal{AK}}_3$-manifolds. It will lead to the proof of the theorem 1.1, hence giving 
a simplified algebraic form of the torsion tensor of an almost K\"ahler structure satisfying the third Gray condition on curvature. \par
Throughout this section $(M^{2n},g,J)$ will be an almost K\"ahler in the class ${\cal{AK}}_3$. An important associated object is the K\"ahler nullity of 
$(g,J)$, the vector subbundle of $TM$ defined by $H=\{v \ \mbox{in} \ TM : \eta_v=0 \}$. 
Note that, a priori, $H$ need not to have constant rank over $M$. However, this is true locally, in the following 
sense. Call a point $m$ of $M$ {\it{regular }}
if the rank of $\eta$ attains a local maximum at $m$. Using standard continuity arguments, it follows that around each regular point, the rank of $\eta$, and 
hence that of $H$ is constant in some open subset. It is also easy to see that the set of regular point is dense 
in $M$, provided that the manifod is connected. As we are concerned with the local (in some neighbourhod of a regular 
point ) structure of ${\cal{AK}}_3$ -manifolds we can assume, without loss of generality, that $H$ has constant rank over $M$. 
This assumption will be made in the whole rest of this paper.\par
We have therefore a $J$-invariant and orthogonal decomposition 
$$ TM={\cal{V}} \oplus H$$
where we define the distribution ${\cal{V}}$ to be the orthogonal of $H$ in $TM$. Let us note that it is straightforward consequence of the definitions (in conjunction with 
the almost K\"ahler condition (2.1)) to check that 
${\cal{V}}=T(TM, TM)$. Then proposition 3.2 can be stated in the equivalent form : 
\begin{nr} \hfill 
$ \overline{\nabla}_V \eta=0\hfill $
\end{nr}
for all $V$ in ${\cal{V}}$. We investigate below some of the most immediate consequences of (4.1). 
\begin{lema}
(i) For all $V,W$ in ${\cal{V}}$ and $X$ we have that $\overline{\nabla}_VW$ belongs to ${\cal{V}}$ and $\overline{\nabla}_VX$ belongs to 
$H$. \\
(ii) The distribution ${\cal{V}}$ is integrable. \\
(iii) $\overline{R}(V_1,V_2).\eta=0$ for all $V_i$ in ${\cal{V}}, i=1,2$.
\end{lema}
{\bf{Proof} } : \\
(i) We have $(\overline{\nabla}_V \eta)(X,U)=0$ for all $U$ in $TM$ or further, since $H$ is the K\"ahler nullity of $(g,J)$, 
$\eta_{\overline{\nabla}_VX}U=0$. This guarantees that $\overline{\nabla}_VX$ belongs to $H$ and since $\overline{\nabla}$  is 
a metric connection it follows that $\overline{\nabla}_VW$ is in ${\cal{V}}$. \\
(ii) follows now by (i) and the (easy to check) fact that $T({\cal{V}}, {\cal{V}}) \subseteq {\cal{V}}$. To get (iii) it suffices to derivate (4.1) in the vertical directions.
$\blacksquare$ \\ \par
Let us set a notational convention, to be used intensively in the present and the next section
and intended to improve presentation. If $E$ and $F$ and vector subbundles of $TM$ and $Q$ is a 
tensor of type $(2,1)$, we will denote by $Q(E,F)$ (or $Q_{E}F$) the subbundle of $TM$ generated 
by elements of the form $Q(u,v)$ where $u$ belongs to $E$ and $v$ is in $F$. We will also 
denote by $<E,F>$ the product of two generic elements of $E$ and $F$ respectively. \par 
The starting point of our discussion will be the following result which was one of the key ingredients leading to the algebraic 
reduction of the torsion for manifolds in the class ${\cal{AK}}_2$. 
\begin{pro} We have an orthogonal, $J$-invariant and $\overline{\nabla}$-parallel (inside ${\cal{V}}$) decomposition 
${\cal{V}}={\cal{V}}_0 \oplus {\cal{V}}_1$ which further satisfies 
\begin{nr} \hfill
$ T({\cal{V}}_0, {\cal{V}}_0)={\cal{V}}_0, \ \eta_{{\cal{V}}_1}{\cal{V}}_1 \subseteq H \  \mbox{and} \ \eta_{{\cal{V}}_1}{\cal{V}}_0=0. \hfill$
\end{nr}
\end{pro}
{\bf{Proof}} : \\
The proof is exactly the same as that given in \cite{Nagy} in the context of ${\cal{AK}}_2$, as the last mentioned proof uses in fact 
only conditions shared by the ${\cal{AK}}_3$ class. In fact, the result holds even more generally, namely for almost K\"ahler manifolds satisfying 
(4.1). 
$\blacksquare$ \\ \par
Since ${\cal{AK}}_3$-geometry is much less rigid than the ${\cal{AK}}_2$-one, we need a finer analysis to produce the proof of the theorem 1.1.
Our discussion needs several preliminary lemmas. First of all, let we introduce the {\it{configuration}} tensor $A : H \times H \to {\cal{V}}$ by 
setting 
$$ \overline{\nabla}_XY=\tilde{\nabla}_XY+A_XY$$ 
for all $X,Y$ in $H$. It is clear that $[A_X,J]=0$ whenever $X$ is in $H$. We also define $B : H \times {\cal{V}} \to H$ by 
$$ \overline{\nabla}_XV=\tilde{\nabla}_XV+B_XV$$ 
for all $X,V$ in $H$ and ${\cal{V}}$ respectively. Obviously, $<A_XY, V>=-<B_XV,Y>$.
\begin{lema}
The distribution ${\cal{V}}_0$ is $\overline{\nabla}$-parallel. 
\end{lema}
{\bf{Proof}} : \\
Let us denote by $A^{+}$ and $A^{-}$ the symmetric resp. skew-symmetric components of $A$. If $V,W$ and $X,Y$ are in ${\cal{V}}$ and $H$ 
respectively we must have (cf. (3.1)) : 
$$\biggl [ (d_{\overline{\nabla}}u)(X,V) \biggr ] (W,Y)=\biggl [ (d_{\overline{\nabla}}u)(W,Y) \biggr ] (X,V).$$
Taking into account the parallelism of $\eta$ over ${\cal{V}}$ (i.e. (4.1)) it follows that 
$$<(\overline{\nabla}_X\eta)(V,W),Y>=-<(\overline{\nabla}_Y\eta)(W,X),V>=
<(\overline{\nabla}_Y\eta)(W,V),X>.$$
We antisymmetrise in $V,W$ and we arrive at 
$$<(\overline{\nabla}_XT)(V,W),Y>+<(\overline{\nabla}_Y T)(V,W),X>=0.$$ 
Since $T(V,W)$ belongs to ${\cal{V}}$ this is clearly equivalent with 
$<T_VW, A_X^{+}Y>=0$ so as to obtain that $A_X^{+}Y$ is in ${\cal{V}}_1$. On the other side, we know that (cf. (3.1))
$$(d_{\overline{\nabla}}u)(X,Y)(V,W)=(d_{\overline{\nabla}}u)(V,W)(X,Y)$$ 
and once again the parallelism of the torsion over ${\cal{V}}$ leads us to 
$<(\overline{\nabla}_X\eta)(Y,V)-(\overline{\nabla}_Y\eta)(X,V), W>=0$. But this is to say that $<\eta_{A^{-}_XY}V,W>=0$ and the 
almost K\"ahler condition ensures that $<A_X^{-}Y, T_VW>=0$ hence $A_X^{-}Y$ belongs to ${\cal{V}}_1$. We showed that 
$A_XY$ belongs to ${\cal{V}}_1$ and this implies that $B_X$ vanishes on ${\cal{V}}_0$ for all $X$ in $H$. It follows that 
$\overline{\nabla}_XV$ is in ${\cal{V}}$ for all $X$ in $H$ and $V$ in ${\cal{V}}_0$. Therefore $(\overline{\nabla}_XT)(V_0,W_0)$ belongs to 
${\cal{V}}$ for all $V_0,W_0$ in ${\cal{V}}_0$, but since $<\overline{\nabla}_X \eta)(V,W),U>=0$ for all $V,W,U$ in ${\cal{V}}$ (again by (3.1) and (4.1))
we see that $(\overline{\nabla}_XT)(V_0,W_0)$ is horizontal, hence equal to 
zero and the $\overline{\nabla}$-parallelism of ${\cal{V}}_0$ is now immediate.
$\blacksquare$  \\ \par
Getting closer to the geometric link between ${\cal{V}}_0$ and its orthogonal complement in $TM$ requires some curvature informations, and these 
are provided by the second Bianchi identity. 
\begin{lema}
Let $V_i, 1 \le i \le 3$ be in ${\cal{V}}$ and $X$ in $H$. We have : \\
(i) $\overline{R}(V_1, V_2, V_3, X)=0$ \\
(ii) $\overline{R}(X,V_1,V_2, V_3)=-<[\eta_{V_2}, \eta_{V_3}]X, V_1>$ \\
(iii) $(\overline{\nabla}_V \overline{R})(X,V_1, V_2,V_3)=0$.
\end{lema}
{\bf{Proof}} : \\
(i) follows directly from lemma 4.1, (i) and the integrability of ${\cal{V}}$. To obtain (ii) one uses the symmetry 
property of corollary 3.1, (i). Finally, (iii) follows by derivating 
(ii) and taking into account that $\overline{\nabla}_VX$ belongs to $H$ for all $X$ in $H$ and $V$ in ${\cal{V}}$ and (4.1).
$\blacksquare$ \\ \par
We will now use the second Bianchi identity for the canonical Hermitian connection in order to get 
more information about the algebraic properties of $\eta$ with respect to the decomposition $TM={\cal{V}} \oplus H$. 
\begin{pro}
Let $X, V_i, 1 \le i \le 4$ be vector fields on $H$ and ${\cal{V}}$ respectively. We have : \\
(i) 
\begin{nr} \hfill 
$ \overline{R}(\eta_{V_2}X, V_1, V_3, V_4)-
\overline{R}(\eta_{V_1}X, V_2, V_3, V_4)=-<[\eta_{V_3}, \eta_{V_4}]X, T_{V_1}V_2>.\hfill $
\end{nr}
(ii) $(\overline{\nabla}_X \overline{R})(V_1, V_2, V_3, V_4)=0$.
\end{pro}
{\bf{Proof}} : \\
Using the second Bianchi identity we obtain 
$$ \begin{array}{lr}
(\overline{\nabla}_X \overline{R})(V_1, V_2, V_3, V_4)+(\overline{\nabla}_{V_1} \overline{R})(V_2, X,V_3, V_4)+
(\overline{\nabla}_{V_2} \overline{R})(X, V_1, V_3, V_4)+\vspace{2mm} \\
\overline{R}(T_XV_1, V_2, V_3.V_4)+\overline{R}(T_{V_1}V_2, X, V_3,V_4)+
\overline{R}(T_{V_2}X, V_1, V_3,V_4)=0.
\end{array} $$ 
Now, the second and the third terms of this equation are vanishing by lemma 4.3, (iii). It is easy to see that the first 
term is $J$-invariant in $V_1$ and $V_2$ and that all the remaining terms are $J$-anti-invariant in $V_1$ and $V_2$. 
Therefore, (ii) is proven and we obtain : $\overline{R}(T_XV_1, V_2, V_3, V_4)+\overline{R}(T_{V_1}V_2, X, V_3,V_4)+
\overline{R}(T_{V_2}X, V_1, V_3,V_4)=0.$ Since $T({\cal{V}}, {\cal{V}}) \subseteq {\cal{V}}$ it suffices now to use 
lemma 4.3, (ii) to conclude. $\blacksquare$ \\ \par
\begin{lema}
Let $(M^{2n},g,J)$ be an ${\cal{AK}}_3$-manifold admitting an orthogonal and $J$-invariant decomposition $TM=D_1 \oplus D_2$ which is also 
$\overline{\nabla}$-parallel. Then the following hold : \\
(i) $\eta_{D_1}T(D_2,D_2)=0$. \\
(ii) 
$$ \overline{R}(V,W,X,Y)=-<T_{T_WX}V+T_{T_XV}W, Y>.$$
whenever $V,W$ are in $D_1$ and $X,Y$ in $D_2$ respectively.
\end{lema}
{\bf{Proof}} : \\
We will prove both assertions in the same time. Let $V,W$ be in $D_1$ and $X,Y$ belong to $TM$ and $D_2$ respectively. Then using the first Bianchi identity for 
$\overline{\nabla}$ (see lemma 3.2) and the parallelism of $D_i, i=1,2$ we get : 
$$ \overline{R}(V,W,X,Y)+\overline{R}(W,X,V,Y)+\overline{R}(X,V,W,Y)=-<T_{T_VW}X+T_{T_WX}V+T_{T_XV}W, Y>.$$
By the parallelism of $D_1$ the last two terms of the left hand side vanish leading to 
$$ \overline{R}(V,W,X,Y)=-<T_{T_VW}X+T_{T_WX}V+T_{T_XV}W, Y>.$$
Now, the right hand side of this equation and the last two terms its the left hand side are $J$-invariant in $X$ and $Y$ 
while the rest is $J$-anti-invariant in the same variables. We conclude that $<T_{T_VW}X,Y>=0$ and since $X$ is arbitrary 
chosen in $TM$ the use of the almost K\"ahler condition (2.1) leads to the proof of (i). Now (ii) follows from (i) by a simple computation. 
$\blacksquare$ 
\begin{coro} We have $\eta_{{\cal{V}}_1}H={\cal{V}}_1$ and $\eta_{{\cal{V}}_0}{\cal{V}}_1=0$. 
\end{coro}
{\bf{Proof}} : \\
With the facts we'we already discussed, the proof is completely analogous to that of the corresponding result in the 
${\cal{AK}}_2$-setting (see \cite{Nagy}). However, we reproduce it below for the convenience of the reader. Let us prove the first 
assertion. \par
Since $\eta_{{\cal{V}}_1}{\cal{V}}_0=0$ we have that 
$\eta_{{\cal{V}}_1}H \subseteq {\cal{V}}_1$. Consider the decomposition ${\cal{V}}_1=E \oplus F$ with 
$\eta_{{\cal{V}}_1}H=E$ and $F$ the orthogonal complement of $E$ in ${\cal{V}}_1$. From the definition of 
$F$ it follows that $\eta_{{\cal{V}}_1}F$ is orthogonal to $H$ and hence it vanishes (recall 
that $\eta_{{\cal{V}}_1} {\cal{V}}_1 \subseteq H$). Since $T({\cal{V}}_1, {\cal{V}}_1)=0$ it also follows that 
$\eta_F {\cal{V}}_1=0$. This implies that $\eta_F H$, which a subspace of ${\cal{V}}_1$, is orthogonal to 
${\cal{V}}_1$, and hence $\eta_FH=0$. Finally since $\eta_F {\cal{V}}_0=0$ ($F$ lies in ${\cal{V}}_1$) we get 
that $F$ is contained in the K\"ahler nullity of $(g,J)$ and then, of course, $F$=0. \par
To prove the second assertion, we note that, by lemma 4.2 the orthogonal and $J$-invariant decomposition $TM={\cal{V}}_0 \oplus ({\cal{V}}_1 \oplus H)$
is $\overline{\nabla}$-parallel. Then by lemma 4.4, (i) we find 
$$ \eta_{{\cal{V}}_0} T({\cal{V}}_1, H)= \eta_{{\cal{V}}_0} \eta_{{\cal{V}}_1}H=0.$$ 
Combining this with $\eta_{{\cal{V}}_1}H={\cal{V}}_1$ finishes the proof of the lemma.
$\blacksquare$ \\ \par
We need now to state a partial version of lemma 4.4. The proof will be omitted since analogous to that of the previously mentioned lemma.
\begin{lema} Let ${\cal{V}}_0=D_1 \oplus D_2$ be a orthogonal and $J$-invariant decomposition which is also $\overline{\nabla}$-parallel (inside 
${\cal{V}}$). Then $\hat{\eta}_{D_1} T(D_2, D_2)=0$, where by $\hat{\eta}$ we denote the ${\cal{V}}_0$ component of the tensor $\eta$. 
\end{lema}
With this preparatives in mind, we can make a decisive step to the proof of theorem 1.1, provided we notice first the following fundamental fact.
We consider the decomposition 
$${\cal{V}}_0=E_1 \oplus E_2$$ 
where $E_2=\eta_{{\cal{V}}_0}H$ and $E_1$ is the orthogonal complement of $E_2$ in ${\cal{V}}_0$. It is clear that this is a $J$-invariant 
decomposition. 
\begin{lema}We have that $\eta_{E_2}E_1=\eta_{E_1}E_2=0$. Moreover, $T(E_i,E_i)=E_i, i=1,2$.
\end{lema}
{\bf{Proof}} : \\
Using the definition of $E_1$ we get that $\eta_{{\cal{V}}_0}E_1$ is orthogonal to $H$ and since $T(E_1, 
{\cal{V}}_0)$ is contained in ${\cal{V}}_0$ it follows that $\eta_{E_1}{\cal{V}}_0$ is orthogonal to $H$. It follows that 
$\eta_{E_1}H$ is orthogonal to ${\cal{V}}_0$ and using that $\eta_{{\cal{V}}_0}{\cal{V}}_1=0$ we arrive at 
\begin{nr} \hfill 
$ \eta_{E_1}H=0.\hfill $
\end{nr}
The definition of $E_2$ implies then $\eta_{E_2}H=E_2$. Also from their definition and the 
parallelism of $\eta$ in vertical directions it is easy to see that 
$E_1$ and $E_2$ are $\overline{\nabla}$-parallel inside ${\cal{V}}$. Therefore, taking $V_3$ in $E_1$ and 
$V_4$ in $E_2$ in (4.3) we find that 
$$ <[\eta_{V_3}, \eta_{V_4}]X, T_{V_1}V_2>=0$$ 
for all $V_i$ in ${\cal{V}}, i=1,2$ and $X$ in $H$. Taking into account that $T({\cal{V}}_0, {\cal{V}}_0)={\cal{V}}_0$ and 
(4.4) we get that $<\eta_{V_3}\eta_{V_4}X, U>=0$ for all $U$ in ${\cal{V}}_0$. Or $\eta_{E_2}H=E_2$ and then 
$<\eta_{E_1}E_2, {\cal{V}}_0>=0$. But $\eta_{E_1}E_2$ is orthogonal to $H$ by (4.4) and also to ${\cal{V}}_1$ by 
corollary 4.1. We arrive at 
\begin{nr} \hfill
$ \eta_{E_1}E_2=0. \hfill $
\end{nr}
Using (4.5) and the almost K\"ahler condition we get $0=<\eta_{E_1}E_2, E_2>=<T(E_2, E_2),E_1>$ ensuring 
 that $T(E_2, E_2) \subseteq E_2$. Again by (4.5) we obtain 
$0=<\eta_{E_1}E_2, E_1>$ that is $\eta_{E_1}E_1$ is orthogonal to $E_2$ and antisymmetrising we also get 
that $T(E_1, E_1) \subseteq E_1$. Thus $<T(E_1, E_1), E_2>=0$ and through the almost K\"ahler condition we obtain 
that $<\eta_{E_2}E_1, E_1>=0$ meaning that $T(E_1, E_2)=\eta_{E_2}E_1 \subseteq E_2$. Or ${\cal{V}}_0=T(
{\cal{V}}_0, {\cal{V}}_0)$, in other words $E_1 \oplus E_2=T(E_1, E_1) +T(E_1,E_2)+T(E_2,E_2)$. Then the above 
inclusions, followed by a dimension argument yield to $T(E_1,E_1)=E_1$. Now lemma 4.5 gives $\hat{\eta}_{E_2}T(E_1, 
E_1)=0$ or further $\hat{\eta}_{E_2}E_1=0$. But $\eta_{E_2}E_1=T(E_1,E_2)$ is contained in ${\cal{V}}_0$ showing that 
$\eta_{E_2}E_1=0$. This finishes the proof of the lemma.
$\blacksquare$ \\ \par
We let now the curvature tensor get into the play by proving the preparatory : 
\begin{lema}
We have : 
$$( \overline{\nabla}_U \overline{R})(V_1,V_2,V_3,V_4)=0$$
for all $U$ in $TM$ and $V_i, 1 \le i \le 4$ in $E_2$. 
\end{lema}
{\bf{Proof}} : \\
If $U$ is in $H$ the result follows by the $\overline{\nabla}$-parallelism of ${\cal{V}}_0$ and proposition 4.2, (ii). Let us suppose now that $U$ belongs to 
${\cal{V}}$. Deriving (4.3) in the direction of $U$ and taking into account the symmetry property of the corollary 3.1, (i) and the parallelism of the torsion over 
${\cal{V}}$ we find that : 
\begin{nr} \hfill 
$ (\overline{\nabla}_U \overline{R})(V_3,V_4,\eta_{V_2}X, V_1)=
(\overline{\nabla}_U \overline{R})(V_3,V_4,\eta_{V_1}X, V_2) \hfill $
\end{nr}
for all $X$ in $H$ and $V_i$ in ${\cal{V}}, 1 \le i \le 4$. Now, by lemma 4.1, (iii), we know that 
$$\overline{R}(V_3,V_4)(\eta_{V_2}X)=
\eta_{V_2}\overline{R}(V_3,V_4)X+\eta_{\overline{R}(V_3,V_4)V_2}X$$
and furthermore, since ${\cal{V}}_0$ is 
$\overline{\nabla}$-parallel, we find that $\overline{R}(V_3,V_4)X=0$ for all $V_3, V_4$ in ${\cal{V}}_0$, accordingly to 
lemma 4.4, (ii). Therefore 
$$\overline{R}(V_3,V_4)(\eta_{V_2}X)=
\eta_{V_2}\overline{R}(V_3,V_4)X$$
for all $V_2, V_3, V_4$ in ${\cal{V}}_0$. Deriving the last equation in the direction of $U$ and invoking again the parallelism 
of $\eta$ in the vertical direction we find after an easy computation (using at its end the almost K\"ahler condition) : 
$$(\overline{\nabla}_U \overline{R})(V_3,V_4,\eta_{V_2}X, V_1)=-(\overline{\nabla}_U \overline{R})(V_3,V_4,\eta_{V_2}X, 
V_1). $$
Together with (4.6) and the fact that $\eta_{E_2}H=E_2$, this leads easily to the proof of the lemma.
$\blacksquare$ \\ \par
Before proceeding to the proof of theorem 1.1 we note that 
making the statement of theorem 1.1 precise requires  the following : 
\begin{defi}
Let $(M^{2n},g,J)$ belong to the class ${\cal{AK}}_3$. It is said to be special iff $\eta_{{\cal{V}}}{\cal{V}} \subseteq H$.
\end{defi}
Note that any $4$-dimensional almost K\"ahler manifold naturally satisfies the algebraic condition stated below. 
We would also like to point out the difference with the definition of special ${\cal{AK}}_2$-manifolds 
(see \cite{Nagy}) where an easy argument allowed reduction to the more pleasant condition $\eta_{{\cal{V}}}{\cal{V}}=H$. 
In the ${\cal{AK}}_3$-case, this reduction is no more an obvious fact. \par 
Another definition required to explain the statement of theorem 1.1 is : 
\begin{defi}
Let $(M^{2n},g,J)$ be an almost K\"ahler with parallel torsion. It is called strict if $\nabla_vJ=0$ implies $v=0$ whenever 
$v$ belongs to $TM$.
\end{defi}
In other words the K\"ahler nullity of of a strict almost K\"ahler manifold with parallel torsion vanishes identically. We are now in position to 
prove the main result of this paper. \\
$\\$
{\bf{Proof of theorem 1.1}} : \\
We are going to show that $E_2=0$, in other words $\eta_{{\cal{V}}_0}H=0$. To this effect 
we consider the partial Hermitian Ricci curvature tensor $\overline{\rho} : {\cal{V}}_0 \to 
{\cal{V}}_0$ defined by 
$$ \overline{\rho} =\sum \limits_{v_k \in E_2} \overline{R}(v_k, Jv_k)$$
where $\{ v_k \}$ is an arbitrary orthonomal basis in $E_2$. 
Using lemma 4.4, (i) we obtain that 
$$\overline{R}(V,W,X,Y)=0$$
for all $V,W$ in ${\cal{V}}_0$ and $X,Y$ in $H$. Then the parallelism of the torsion over ${\cal{V}}$ leads 
easily to $\overline{\rho}(\eta_VX)=\eta_{\overline{\rho}(V)}X$ for all $V$ in ${\cal{V}}_0$ and $X$ in $H$. We write 
now $\overline{\rho}=SJ$ where $S$ is a symmetric, $J$-invariant endomorphism of ${\cal{V}}_0$. Then the previous 
relation reads : 
\begin{nr} \hfill 
$S(\eta_VX)=-\eta_{SV}X \hfill $
\end{nr}
whenever $V$ is in ${\cal{V}}_0$ and $X$ in $H$. Using lemma 4.7, it is easy to see that the 
restriction of $S$ to $E_2$ is $\overline{\nabla}$-parallell, hence has a globally defined spectral decomposition 
and constant eigenvalues. We are now going to show that $S=0$. By contradiction, let us assume that the restriction 
of $S$ to $E_2$ has a non-zero eigenvalue $\lambda$, and denote the corresponding 
eigenspace by ${\cal{W}}_1$. Set 
${\cal{W}}_2=\eta_{{\cal{W}}_1}H$ and note that by (4.7) we have that ${\cal{W}}_2 \subseteq Ker(S + \lambda)$. Then ${\cal{W}}_1$ and 
${\cal{W}}_2$ are orthogonal, and again by (4.7) it follows that 
$\eta_{{\cal{W}}_2}H \subseteq {\cal{W}}_1.$
Let ${\cal{W}}$ be the orthogonal complement of 
${\cal{W}}$ in ${\cal{V}}_0$. Then by (4.7) we get that $\eta_{{\cal{W}}}H \subseteq {\cal{W}}$. Now 
$\eta_{E_2}H=E_2$, hence a dimension argument 
shows that we must have 
\begin{nr} \hfill 
$\eta_{{\cal{W}}_2}H={\cal{W}}_1$ and $\eta_{{\cal{W}}}H={\cal{W}} . \hfill $ 
\end{nr}
From the $\overline{\nabla}$-parallelism of $S$ it follows immediately that the distributions ${\cal{W}}, i=1,2$ and 
${\cal{W}}$ are equally $\overline{\nabla}$-parallel. With this fact in mind we will now make use of  
proposition 4.2, (i). Taking $V_3$ in ${\cal{W}}$ and $V_4$ in ${\cal{W}}_1$ in the relation (4.3) we obtain 
$$ <[\eta_{V_3}, \eta_{V_4}]X, T(V_1,V_2)>=0$$ 
for all $V_1,V_2$ in ${\cal{V}}$ and $X$ in $H$. Since $T({\cal{V}}_0, {\cal{V}}_0)={\cal{V}}_0$ it follows further 
that $<[\eta_{V_3}, \eta_{V_4}]X, U>=0$ for all $U$ in ${\cal{V}}_0$. Now an invariance 
(with respect to $J$) argument in the variables $V_3, X$ for example yields to 
$<\eta_{V_3}\eta_{V_4}X,U>=<\eta_{V_4}\eta_{V_3}X,U>=0$.  Using (4.8), we get that $\eta_{{\cal{W}}}{\cal{W}}_1$ 
and $\eta_{{\cal{W}}_1}{\cal{W}}$ are 
orthogonal to ${\cal{V}}_0$. But these spaces are orthogonal to ${\cal{V}}_1$ (see corollary 4.1) and also to $H$ by 
(4.8). Therefore 
\begin{nr} \hfill 
$\eta_{{\cal{W}}}{\cal{W}}_1=\eta_{{\cal{W}}_1}{\cal{W}}=0 \hfill $
\end{nr}
and in a completely analogous manner one arrives to 
\begin{nr} \hfill
$ \eta_{{\cal{W}}_2}{\cal{W}}_1=\eta_{{\cal{W}}_1}{\cal{W}}_2=\eta_{{\cal{W}}}{\cal{W}}_2=\eta_{{\cal{W}}_2}{\cal{W}}=0. \hfill $
\end{nr}
Using the almost K\"ahler condition (2.1) it is easy to derive from these properties that 
$T({\cal{W}}, {\cal{W}}) \subseteq {\cal{W}}$ and 
$ T({\cal{W}}_i, {\cal{W}}_i) \subseteq {\cal{W}}_i, i=1,2$ as well as 
$T({\cal{W}}, {\cal{W}}_i)=0, i=1,2$. But $T(E_2, E_2)=E_2$ and 
$E_2={\cal{W}}_1 \oplus {\cal{W}}_2 \oplus {\cal{W}}$ hence a dimension 
argument shows that $T({\cal{W}}_i, {\cal{W}}_i) ={\cal{W}}_i, i=1,2$ and 
$T({\cal{W}}, {\cal{W}}) \subseteq {\cal{W}}$. In particular : 
\begin{nr} \hfill 
$ T({\cal{W}}_1, {\cal{W}}_1)={\cal{W}}_1. \hfill $
\end{nr}
Now, the space $\eta_{{\cal{W}}_1} {\cal{W}}_1 $ is orthogonal to ${\cal{W}}_2 \oplus {\cal{W}} $ by (4.9) and (4.10), 
to $E_1$ by lemma 4.6, to ${\cal{V}}_1$ by corollary 4.1 and finally to $H$ by (4.8). We are lead to : 
\begin{nr} \hfill 
$ \eta_{{\cal{W}}_1} {\cal{W}}_1 \subseteq {\cal{W}}_1.\hfill $ 
\end{nr}
Using again the parallelism of the torsion over ${\cal{V}}$ we find that 
$$ \overline{\rho}(\eta_vw)=\eta_{\overline{\rho}v}w+\eta_{v}\overline{\rho}v $$
for all $v,w$ in ${\cal{W}}_1$. In view of (4.12) and of the fact $\overline{\rho}$ acts on 
${\cal{W}}_1$ as $\lambda J$ we obtain that 
$\eta_{{\cal{W}}_1} {\cal{W}}_1$ vanishes and then ${\cal{W}}_1$ must be 
vanishing too by (4.11). We obtained a contradiction leading to 
the fact that $\overline{\rho}=0$ on $E_2$. But this means that any integral manifold of $E_2$ is an almost K\"ahler manifold with parallel torsion and vanishing Hermitian 
Ricci tensor. As shown in \cite{Nagy} such a manifold has to be K\"ahler hence the 
fact that $T(E_2, E_2)=E_2$ ensures the vanishing of $E_2$. Therefore we have that $\eta_{{\cal{V}}_0}H=0$ and this implies immediately $\eta_{{\cal{V}}_0}{\cal{V}}_0 
\subseteq {\cal{V}}_0$. Since ${\cal{V}}_0$ is a $\overline{\nabla}$-parallel distribution, the last conditions are in fact telling us that ${\cal{V}}_0$ is 
$\nabla$-parallel and some straightforward considerations are now finishing the proof.
$\blacksquare$ \\ 
\par
To prove the corollary 1.1 one notices first that if the torsion is parallel, then the K\"ahler nullity is parallel with respect 
to the first Hermitian connection and hence the decomposition of theorem 1.1 holds globally. Moreover, in 
\cite{Nagy3} it was shown that a special (in the sense of definition 4.1) ${\cal{AK}}_2$-manifold is locally the product 
of a K\"ahler  manifold and a space of type I, hence finishing the proof of the corollary.